\documentclass[12pt]{article}
\usepackage{latexsym, amsmath, amsfonts, amscd, amssymb, verbatim, longtable}
\usepackage{multirow, booktabs}

\usepackage{color}

\def\tto{\;{\lower 1pt \hbox{$\rightarrow$}}\kern -10pt
	\hbox{\raise 2pt \hbox{$\rightarrow$}}\;}

\begin{document}
	\pagestyle{myheadings}
	
	\newtheorem{Theorem}{Theorem}[section]
	\newtheorem{Proposition}[Theorem]{Proposition}
	\newtheorem{Remark}[Theorem]{Remark}
	\newtheorem{Lemma}[Theorem]{Lemma}
	\newtheorem{Corollary}[Theorem]{Corollary}
	\newtheorem{Definition}[Theorem]{Definition}
	\newtheorem{Example}[Theorem]{Example}
	\renewcommand{\theequation}{\thesection.\arabic{equation}}
	\normalsize
	\setcounter{equation}{0}
	\title{\bf Convergence of a Solution Algorithm in Indefinite Quadratic Programming\footnote{The authors were supported respectively by the Faculty of Information Technology of Hanoi University of Industry, the National
				Research Foundation of Korea (NRF) grant funded by the Korea
				government (MEST) No.2015R1A3A2031159, and the National Foundation for Science $\&$ Technology Development (Vietnam). The first and third authors thank the Sungkyunkwan University for supporting their research stays in Suwon.}}
	
	%\author
	\medskip
	\author{T.H. Cuong\footnote{Department of Computer Science, Faculty of Information Technology, Hanoi University of Industry, 298 Cau Dien Road, Bac Tu Liem District, Hanoi, Vietnam; email: tranhungcuong@haui.edu.vn.},\ \, Y. Lim\footnote{Department of Mathematics, Sungkyunkwan University, Suwon 440-746, South Korea; email: ylim@skku.edu},\ \, N.D. Yen\footnote{Institute of Mathematics, Vietnam Academy of
			Science and Technology, 18 Hoang Quoc Viet, Hanoi 10307, Vietnam;
			email: ndyen@math.ac.vn.}}\maketitle
	\date{}%\small\today}

\medskip
\begin{quote}
\noindent {\bf Abstract.}  It is proved that, for an indefinite quadratic programming problem under linear constraints,  any iterative sequence generated by the {\it Proximal DC decomposition algorithm} $R$-linearly  converges to a Karush-Kuhn-Tucker point, provided that the problem has a solution. Another major result of this paper says that DCA sequences generated by the algorithm converge to a locally unique solution of the problem if the initial points are taken from a suitably-chosen neighborhood of it. To deal with the implicitly defined iterative sequences, a local error bound for affine variational inequalities and  novel techniques are used. Numerical results together with an analysis of the influence of the decomposition parameter, as well as a comparison between the Proximal DC decomposition algorithm and the\textit{ Projection DC decomposition algorithm}, are given in this paper. Our results complement a recent and important paper of  Le Thi, Huynh, and Pham Dinh (J. Optim. Theory Appl. {\bf 179} (2018), 103--126).


\end{quote}

%\newpage
\section{Introduction}

 The importance of the indefinite quadratic programming problem under linear constraints (IQP for brevity) in optimization theory and its various applications is well known (see, e.g., \cite{Bomze_1998,Gould_Toint}). 
 
 \medskip
 For the solution existence, structure of the solution set, necessary and sufficient optimality conditions, and stability of this problem, the reader is referred to \cite{LeeTamYen_book} and the references therein.

\medskip
Numerical methods for solving IQP have been addressed in many research works; see, e.g.,  \cite{Bomze_Danninger_1994,Cambini_Sodini_2005,PhamDinh_LeThi_2,PhamDinh_LeThi98,PhamDinh_LeThi_3,PhamDinh_LeThi_Akoa,Ye89,Ye92,Ye97}. Note that most of the known algorithms yield just stationary points (that is, the Karush-Kuhn-Tucker points, or KKT points for short), or local minimizers. In other words, most of the known algorithms are \textit{local solution methods}. Since the IQP is NP-hard (see \cite{Pardalos_Vavasis_1991} and also \cite{Bomze_Danninger_1994}), finding its \textit{global solutions} remains a challenging question.

\medskip
We are interested in studying and implementing two methods to solve the IQP, that are based on a general scheme for solving DC (Difference-of-Convex-functions) programs due to Pham Dinh and Le Thi \cite{PhamDinh_LeThi_AMV97,PhamDinh_LeThi98} (see also \cite{LeThi_PhamDinh_AOR05,PhamDinh_LeThi_4}). A combination of DCA (DC Algorithms) with interior point techniques for solving large-scale nonconvex quadratic programming has been proposed in \cite{PhamDinh_LeThi_Akoa}. The two DC decompositions suggested in \cite{PhamDinh_LeThi_Akoa} are the projection DC decomposition and the proximal DC decomposition. They lead to two algorithms for solving the IQP: the {\it Projection DC decomposition algorithm} (Algorithm A) and the {\it Proximal DC decomposition algorithm} (Algorithm B); see \cite{ATY2,PhamDinh_LeThi_Akoa}, and Section 2 below. It is worthy to stress the following features of these algorithms:

\vskip6pt
- The algorithm descriptions are simple;

- The implementation is easy;

- No line searches are required.

\medskip
Nevertheless, using the DCA theory one can only assert \cite[Theorem~1]{ATY1} that any cluster point of a DCA sequence generated by the above-mentioned algorithms is a KKT point of the IQP. To be sure that such cluster points do exist, one must establish the boundedness of the DCA sequence. In general,  DCA sequences need not be bounded \cite[Example~1]{ATY1}. But there is a Conjecture \cite[p.~489]{ATY1} saying that \textit{if the IQP has global solutions, then every DCA sequence generated by one of the algorithms A and B must be bounded.} Recently, the Conjecture has been solved in the affirmative for the two-dimensional IQP by Tuan \cite{Tuan_JOTA2015}. To solve it in the general case, Tuan \cite{Tuan_JMAA2015} has used a local error bound for affine variational inequalities and several specific properties of the KKT point set of the IQP which were obtained by Luo and Tseng \cite{Luo_Tseng_1992} (see also Tseng \cite{Tseng_1995} and Luo \cite{Luo_2000}). The main result of \cite{Tuan_JMAA2015} is the following theorem: \textit{If the IQP has a nonempty  solution  set, then every DCA sequence generated by Algorithm A converges $R$-linearly to a KKT point.} 

\medskip
The first aim of the present paper is to prove that \textit{any DCA sequence generated by Algorithm B converges $R$-linearly to a KKT point}. Hence, combining this with Theorem 2.1 from  \cite{Tuan_JMAA2015}, we have a complete solution for the Conjecture in \cite[p.~489]{ATY1}. Our result is obtained by applying some arguments of \cite{Tuan_JMAA2015} and a new technique in dealing with implicitly defined DCA sequences.

\medskip
 By \cite[Theorem~3]{ATY1}, we know that DCA sequences generated by the Algorithm~A converge to a locally unique solution of the IQP if the initial points are taken from a suitably-chosen neighborhood of it. In the terminology of \cite{LeongGoh13}, this means that the locally unique solutions of the IQP are \textit{asymptotically stable} with respect to Algorithm A. The open question of \cite[p.~488]{ATY1} can be reformulated as follows: \textit{Is it true that the locally unique solutions of the IQP are asymptotically stable with respect to Algorithm~B?}

\medskip
The second aim of our paper is to use a novel technique to establish the asymptotical stability of the locally unique solutions with respect to Algorithm~B under a mild additional assumption on the DCA decomposition parameter. It is still unclear to us whether that assumption can be dropped, or not.

\medskip
The third aim of this paper is to analyze the influence of the decomposition parameter on the rates of convergence of DCA sequences and compare the performances of the algorithms A and B upon randomly generated data sets. Numerous numerical tests lead us to the following observations:

\vskip6pt
- For both the the algorithms A and B, the closer is the positive decomposition parameter to the lower bound of the admissible parameter interval, the higher is the convergence rate of DCA sequences;

- Algorithm B is more efficient and more stable than Algorithm A.

\medskip
Our results complement a recent paper of  Le Thi, Huynh, and Pham Dinh \cite{LeThi_Ngai_PhamDinh_JOTA2018}, where by original proofs the authors have obtained a series of important convergence theorems for DCA algorithms, which solve optimization problems with subanalytic data. To be more precise, from Theorems 3.4, 3.5, and 4.2 of \cite{LeThi_Ngai_PhamDinh_JOTA2018} it follows that any DCA sequence generated by Algorithm B converges $R$-linearly to a KKT point, \textit{if the sequence is bounded}. Since the boundedness of DCA sequences cannot be obtained by the Lojasiewicz inequality (see \cite[Theorem~2.1]{LeThi_Ngai_PhamDinh_JOTA2018}) and the related results on Kurdyka-Lojasiewicz properties (see \cite{An_Nam_OPTIM2017} and the references therein), Theorem \ref{Main Result} and its proof are new contributions to the analysis of the existing solution algorithms in indefinite quadratic programming.

\medskip
The interested reader is referred to the comprehensive survey paper of Le Thi and Pham Dinh \cite{LeThi_PhamDinh_MP2018} on the thirty years (1985--2015) of the development of the DC programming and DC algorithms, where as many as 343 research works have been commented and the following remarks have been given: \textit{``DC programming and DCA were the subject of several hundred articles in the high ranked scientific journals and the high-level international conferences, as well as various international research projects, and were the methodological basis of more than 50 PhD theses. About 100 invited symposia/sessions dedicated to DC programming and DCA were presented in many international conferences. The ever-growing number of works using DC programming and DCA proves their power and their key role in nonconvex programming/global optimization and many areas of applications."}

\medskip
The remainder of the paper consists of four sections.  Section 2 describes the DC algorithms of \cite{PhamDinh_LeThi_Akoa}. In Section 3, we study the $R$-linear convergence rate of DCA sequences generated by Algorithm~B. Section 4 establishes a theorem on the asymptotical stability of the locally unique solutions with respect to Algorithm B.  The influence of the decomposition parameter on the rates of convergence of DCA sequences and the performances of the algorithms A and B upon randomly generated data sets are discussed in Section 5. 

\section{Preliminaries}
\setcounter{equation}{0}

Consider \textit{the indefinite quadratic programming problem under linear constraints} (called the IQP in the preceding section):
\begin{eqnarray}\label{QP problem}\min\Big\{f(x):=\frac{1}{2}x^TQx+q^Tx\, :\, Ax\geq b\Big\}, \end{eqnarray}
where $Q\in\mathbb R^{n\times n}$ and $A\in\mathbb  R^{m\times n}$ are given matrices, $Q$  is symmetric, $q\in\mathbb  R^n$ and $b\in\mathbb  R^m$ are arbitrarily given vectors. The constraint set of the problem is  $$C:=\big\{x\in\mathbb R^n\, :\, Ax\geq b\big\}.$$ Since $x^TQx$ is an indefinite quadratic form, the objective function $f(x)$ may be nonconvex; hence (\ref{QP problem}) is a nonconvex optimization problem. 

\medskip Now we describe some standard notations that will be used later on. The unit matrix in $\mathbb R^{n\times n}$ is denoted by $I$. The \textit{eigenvalues} of a symmetric matrix $M\in\mathbb R^{n\times n}$ are ordered in the sequence $\lambda_1(M)\leq ...\leq\lambda_n(M)$ with counting multiplicities. For an index set $\alpha\subset \{1,\dots,m\}$, by $A_\alpha$ we denote the matrix composed by the rows $A_i$, $i\in\alpha$, of $A$. Similarly, $b_\alpha$ is the vector composed by the components $b_i$, $i\in\alpha$, of $b$. The {\it pseudo-face} of $C$ corresponding to $\alpha$ is the set $$\big\{x\in\mathbb R^n\,:\, A_\alpha x=b_\alpha,\ A_{\bar\alpha}x>b_{\bar\alpha}\big\},$$ where $\bar\alpha:=\{1,\dots,m\}\setminus\alpha$. Let $B(x,\varepsilon)$ (resp., $\overline{B}(x,\varepsilon)$) denote the open (resp., closed) ball with center $x$ and radius $\varepsilon>0$. Given $s$ vectors $v^1,\dots,v^s$ in $\mathbb R^n$, we denote by ${\rm pos}\{v^1,\dots,v^s\}$ the closed convex cone generated by $v^1,\dots,v^s$, that is $${\rm pos}\{v^1,\dots,v^s\}= \Big\{v=\sum_{i=1}^s\lambda_iv^i\,:\, \lambda_i\geq 0\ \,{\rm for}\ \,i=1,\dots,s\Big\}.$$ Symbol $\Omega^\perp$ stands for the linear subspace consiting of the vectors which are orthogonal to every vector in $\Omega$. The metric projection of $u\in\mathbb R^n$ onto $C$ is denoted by $P_C(u)$, that is  $P_C(u)\in C$ and $$\big\|u-P_C(u)\big\|=\displaystyle\min_{x\in C} \|u-x\|.$$ The tangent cone to $C$ at $x\in C$ is denoted by $T_C(x)$, i.e., $$T_C(x)=\{t(y-x)\,:\,t\geq 0,\ y\in C\}=\{v\in\mathbb R^n\,:\,A_\alpha v\geq 0\},$$ where $\alpha=\{i\,:\, A_ix=b_i\}$. The normal cone to $C$ at $x\in C$ is denoted by $N_C(x)$, that is 
\begin{eqnarray*} N_C(x)=\big(T_C(x)\big)^* &=& \{\xi\in\mathbb R^n\,:\,\langle\xi,v\rangle\leq 0\ \, \forall v\in T_C(x)\}\\ &=& -{\rm pos}\{A_i\,:\, i\in\alpha\}.
\end{eqnarray*}  

Following \cite{PhamDinh_LeThi_Akoa}, to solve the IQP via a sequence of strongly convex quadratic programs, one decomposes $f(x)$ into the difference of two convex linear-quadratic functions
\begin{eqnarray}\label{DC decomposition of f}f(x)=\varphi(x)-\psi(x)\end{eqnarray}  with $\varphi(x)=\frac{1}{2}x^TQ_1x+q^Tx$ and $\psi(x)=\frac{1}{2}x^TQ_2x$, where $Q=Q_{1}-Q_{2}$, $Q_1$ is a symmetric positive definite matrix and $Q_2$ is a symmetric positive semidefinite matrix. Then (\ref{QP problem}) is equivalent to the \textit{DC program}
$$\min\big\{g(x)-h(x)\,:\,x\in\mathbb R^n\big\}$$ with $g(x):=\varphi(x)+\delta_C(x)$, $h(x):=\psi(x)$, where $\delta_C(x)=0$ for $x\in C$ and $\delta_C(x)=+\infty$ for $x\notin C$ is the indicator function of $C$. Let
$x^0\in\mathbb R^n$ be a given initial point. In accordance with the general solution method of \cite{PhamDinh_LeThi98, PhamDinh_LeThi_4}, at every step $k\geq 0$ one computes $y^k=\big(\nabla h(x^k)\big)^T=Q_2x^k$ and finds the unique solution, denoted by $x^{k+1}$ of the convex minimization problem
\begin{eqnarray*}\label{auxiliary problem Pk}
	\min\Big\{g(x)-[h(x^k)+\langle x-x^k,y^k\rangle]\,:\,x\in\mathbb R^n\Big\}.
\end{eqnarray*} The latter is equivalent to the \textit{strongly convex quadratic program}
\begin{eqnarray}\label{auxiliary problem Pk_2}
\min\Big\{\frac{1}{2}x^TQ_1x+q^Tx-x^TQ_2x^k\,:\,x\in C\Big\}.
\end{eqnarray} 
The obtained sequence $\{x^k\}$ is called {\it the DCA sequence generated by the DC algorithm and the initial point $x^0$.}

\begin{Definition} {\rm For $x\in\mathbb R^n$, if there exists a multiplier $\lambda\in\mathbb R^m$  such that
		\begin{eqnarray*}
			\begin{cases}Qx+q-A^T\lambda=0,\\ Ax\geq b,\ \; \lambda\geq 0,\ \; \lambda^T(Ax-b)=0,\end{cases}\end{eqnarray*} then $x$ is said to be a {\it Karush-Kuhn-Tucker point} (a {\it KKT point}) of the IQP.
	} \end{Definition}
	
	This definition can be rephrased (see, e.g., \cite{LeeTamYen_book}) as follows: If $x\in C$ and 
	\begin{eqnarray}\label{KKT_point_tangent_cone}\langle \nabla f(x),v\rangle=(Qx+q)^Tv\geq 0\quad \forall v\in T_C(x),
	\end{eqnarray} then $x$ is said to be a KKT point of (\ref{QP problem}). Since condition \eqref{KKT_point_tangent_cone} is equivalent to $\langle \nabla f(x),y-x\rangle\geq 0$ for all $y\in C$, $x\in C$ is a KKT point of the IQP in \eqref{QP problem} if and only if it is a solution of the \textit{affine variational inequality}
		\begin{eqnarray}x\in C,\quad \label{AVI}\langle Qx+q,u-x\rangle\geq 0\ \; \forall u\in C.
	\end{eqnarray}
	Denote the KKT point set (resp., the global solution set) of IQP by $C^*$ (resp., $\mathcal{S}$). It is well known (see, e.g., \cite{LeeTamYen_book}) that $\mathcal{S}\subset C^*$. 
	
	\medskip We now recall some basic properties of DCA sequences.
	
	\begin{Theorem}\label{basic properties of DCA sequences} {\rm (See \cite[Theorem 3]{PhamDinh_LeThi_4} and \cite[Theorem 2.1]{PhamDinh_LeThi_Akoa})} Every DCA sequence $\{x^k\}$ generated by the above DC algorithm and an  initial point $x^0\in\mathbb R^n$ has the following properties:
		\begin{itemize}
			\item[{\rm (i)}] $f(x^{k+1})\leq f(x^k)-\displaystyle\frac{1}{2}[\lambda_1(Q_1)+\lambda_1(Q_2)]\|x^{k+1}-x^k\|^2$ for every $k\geq 1$;
			\item[{\rm (ii)}] $\{f(x^k)\}$ converges to an upper bound $f_*$ for the optimal value of {\rm (\ref{QP problem})};
			\item[{\rm (iii)}] Every cluster point $x^*$ of $\{x^k\}$ is a KKT point of {\rm (\ref{QP problem})};
			\item[{\rm (iv)}] If $\displaystyle\inf_{x\in C}f(x)>-\infty$, then $\displaystyle\lim_{k\to \infty}\|x^{k+1}-x^k\|=0$.
		\end{itemize}
	\end{Theorem}
	
	\begin{Remark}\label{Remark 1} {\rm By \cite[Theorem 3]{PhamDinh_LeThi_4}, if $x^0\in C$ then we have the inequality in (i) for every $k\geq 0$. To see this, it suffices to note that $x^0\in C={\rm dom}g:=\{x\,:\, g(x)<+\infty\},$ where $g=\varphi+\delta_C$.}
	\end{Remark}
	
As the smallest eigenvalue $\lambda_1(Q)$ and the largest eigenvalue $\lambda_n(Q)$ of $Q=Q_{1}-Q_{2}$ can be computed easily by some algorithm (for instance, by the Newton-Raphson algorithm in \cite{Stoer_Bulirsch_1980}) or software, next realizations of the DC decomposition~(\ref{DC decomposition of f}) can be done: 
	
	(a) $Q_1:=\rho I$, $Q_2:=\rho I-Q$, where $\rho$ is a positive real value satisfying the condition $\rho\geq\lambda_n(Q)$;
	
    (b) $Q_1:=Q+\rho I$, $Q_2:=\rho I$, where $\rho$ is a positive real value satisfying the condition $\rho>-\lambda_1(Q)$.
    
    \medskip The number $\rho$ is called \textit{the decomposition parameter}. The following algorithms appear on the basis of (a)
and (b), respectively.
	
	\medskip
	{\bf Algorithm A.} (\textit{Projection DC decomposition algorithm}) Fix a positive number $\rho\geq\lambda_n(Q)$ and choose an initial point $x^0\in\mathbb R^n$. For every $k\geq 0$, compute the point \begin{eqnarray*}\label{Iteration formula of Algorithm A} x^{k+1}:=P_C\Big(x^k-\frac{1}{\rho}(Qx^k+q)\Big)\end{eqnarray*} which is the unique solution of (\ref{auxiliary problem Pk_2}),
	where $Q_1=\rho I$ and $Q_2:=\rho I-Q$. The latter can be rewritten in the form 
	$$\min\Big\{\Big\|x-\frac{1}{\rho}(y^k-q)\Big\|^2\, :\, Ax\geq b\Big\}$$ with $y^k:=(\rho I-Q)x^k$.
	
	\medskip
	{\bf Algorithm B.}  (\textit{Proximal DC decomposition algorithm}) Fix a positive number $\rho> -\lambda_1(Q)$ and choose an initial point $x^0\in\mathbb R^n$. For any $k\geq 0$, compute the unique solution, denoted by the point $x^{k+1}$, of the strongly convex quadratic minimization problem
	\begin{eqnarray}\label{Minimization problem of Algorithm B}\min\Big\{\psi(x):=\frac{1}{2}x^TQx+q^Tx+\frac{\rho}{2}\|x-x^k\|^2\, :\, Ax\geq b\Big\}.\end{eqnarray} (Note that, up to adding a real constant, the objective function of (\ref{Minimization problem of Algorithm B}) can be written as $\frac{1}{2}x^TQ_1x+q^Tx-x^TQ_2x^k$, where $Q_1=Q+\rho I$ and $Q_2=\rho I$.)
	
	\medskip
	Let $\{x^k\}$ be a DCA sequence generated by one of the last two algorithms and an initial point $x^0$. If $\{x^k\}$ is bounded, then it has a convergent subsequence $x^{k_j}\to x^*$. According to Theorem \ref{basic properties of DCA sequences}, $x^*$ is a KKT point of IQP. Since one wants to find a global solution, one has to {\it restart} the algorithm if  $x^*\notin\mathcal{S}$. To do so, we must find some $u\in C$ such that $f(u)<f(x^*)$, put $x^0=u$ and construct a new DCA sequence. If the latter is again bounded, one finds a new KKT point $\bar u\in C^*$ with $f(\bar u)\leq f(u)<f(x^*)$ (see Theorem \ref{basic properties of DCA sequences}). The process is continued until finding a point $x_*\in\mathcal{S}$. Since the distinct values of $f$ on $ C^*$ does not exceed $2^m$ (see \cite[Lemma 4]{Bomze_Danninger_1994}), the upper bound for the number of restarts of any DC algorithm is $2^m$. 

\section{Convergence Theorem}
	
	As noted in Section 2, the KKT point set $C^*$ of \eqref{QP problem}  is the solution set of the affine
	variational inequality \eqref{AVI}, so $C^*$ is the union of
	finitely many polyhedral convex sets (see, e.g.,
	\cite[Lemma~3.1]{Luo_Tseng_1992} and \cite [Sections 3.1 and 5.3]
	{LeeTamYen_book}). In particular, $C^*$ has finitely many connected
	components. Since the solution set of \eqref{QP problem} is a subset of $C^*$, if the former is nonempty then  $C^*\ne \emptyset$. 	For any given subset $M \subset \mathbb{R}^n$, by
	$d(x,M):=\inf\{\|x-y\|\,:\,y\in M\}$ one denotes the distance from
	$x\in \mathbb{R}^n$.
	
	\medskip
				We will need two lemmas for our purposes. Next lemma gives a local error bound for the distance form a feasible point $x\in C$ to $C^*$.
		
		\begin{Lemma}\label{Lemma 1} {\rm (\cite[Lemma 2.1]{Tuan_JMAA2015}; cf. \cite[Lemma 3.1]{Luo_Tseng_1992})}\ \,
			For any $\rho>0$, if $C^*\neq\emptyset$, then there exist scalars $\varepsilon >0$ and $\ell>0$ such that
			\begin{equation}\label{1new}
			\begin{array}{rl}
			d(x,C^*)\le \ell \Big\|x-P_C\Big(x-\dfrac{1}{\rho}(Qx+q)\Big)\Big\|
			\end{array}
			\end{equation}
			for all $x\in C$ with
			\begin{equation}\label{3new}
			\begin{array}{rl}
			\Big\|x-P_C\Big(x-\dfrac{1}{\rho}(Qx+q)\Big)\Big\|\le \varepsilon.
			\end{array}
			\end{equation}
		\end{Lemma}
			
		\begin{Lemma}\label{Lemma 3} {\rm (\cite[Lemma 3.1]{Luo_Tseng_1992}; see also \cite[Lemma 2.2]{Tuan_JMAA2015})}\ \,
			Let $C_1, C_2, \cdots, C_r$ denote the connected components of $C^*$. Then we have
			$$ C^*=\displaystyle\bigcup_{i=1}^r C_i,$$
			and the following properties are valid:
			
			{\rm(a)} each $C_i$ is the union of finitely many polyhedral convex sets;
			
			{\rm(b)} the sets $C_i$, $i=1,\ldots r$, are properly separated each from others, that is, there exists $\delta>0$ such that if  $i\ne j$ then
			$$d(x,C_j)\geq \delta\quad \forall x\in C_i;$$
			
			{\rm(c)} $f$ is constant on each $C_i$.
		\end{Lemma}

	The necessary and sufficient condition for $x^{k+1}$ to be the unique solution of \eqref{Minimization problem of Algorithm B} is the following
	$$\langle \nabla\psi(x^{k+1}), x-x^{k+1}\rangle\geq 0\quad \forall x\in C,$$
	where $\nabla\psi(x^{k+1})=Qx^{k+1}+q+\rho x^{k+1}-\rho x^k$. Equivalently, $x^{k+1}$ is the unique solution of the strongly monotone affine variational inequality given by the affine operator $x\mapsto (Q+\rho I)x+q-\rho x^k$ and the polyhedral convex set $C$. Therefore, applying Theorem 2.3 from \cite[p.~9]{KS1980} we see that $x^{k+1}$ is \textit{the unique fixed point} of the map
	$G_k(x):=P_C(x-\mu(Mx+q^k))$, where $\mu>0$ is arbitrarily chosen,  $M:=Q+\rho I$, and $q^k:=q-\rho x^k$. In what follows, we choose $\mu=\rho^{-1}$. Then
	\begin{eqnarray}\label{Iteration formula of Algorithm B} x^{k+1}=P_C\Big(x^{k+1}-\frac{1}{\rho}(Mx^{k+1}+q^k)\Big)\end{eqnarray}
	
	The convergence and the rate of convergence of 	Algorithm  B, the Proximal DC decomposition algorithm, can be formulated as follows.
	
	\begin{Theorem}\label{Main Result} If \eqref{QP problem} has a solution,
		then for each $x^0\in \mathbb{R}^n$, the DCA sequence $\{x^k\}$ constructed by
		Algorithm  B converges
		$R$-linearly to a KKT point of \eqref{QP problem}, that is, there
		exists $x^*\in C^*$ such that
		$$\limsup_{k\to \infty}\|x^k-x^*\|^{1/k}<1.$$
	\end{Theorem}
	{\textit{Proof} Since \eqref{QP problem} has a solution, $C^*\neq\emptyset$. Hence, by Lemma \ref{Lemma 1} there exist $\ell>0$ and $\varepsilon>0$ such that \eqref{1new} is fulfilled for any $x$ satisfying \eqref{3new}. As $\displaystyle\inf_{x\in C}f(x)>-\infty$, assertion~(iv) of Theorem \ref{basic properties of DCA sequences} gives
	\begin{equation}\label{assym_stab}\lim_{k\to\infty}\|x^{k+1}-x^k\|=0.	\end{equation} Choose $k_0\in \mathbb{N}$ as large as $\|x^{k+1}-x^k\|<\varepsilon$ for all $k\ge k_0$. 
	
	If it holds that
		\begin{equation}\label{ineq1n}
		\|x^{k+1}-P_C(x^{k+1}-\frac{1}{\rho}(Qx^{k+1}+q))\|\leq\varepsilon\quad
		\forall k\geq k_0,
		\end{equation} 
			then by \eqref{1new} one has 
			\begin{equation} \label{ineq1}
			d(x^{k+1},C^*)\le \ell \|x^{k+1}-P_C\Big(x^{k+1}-\frac{1}{\rho}(Qx^{k+1}+q)\Big)\|\quad\forall k\geq k_0.
			\end{equation}
		To obtain \eqref{ineq1n}, for any $k\geq k_0$, we recall that 
		\begin{equation}\label{fixed_point}x^{k+1}=G_k(x^{k+1})=P_C\Big(x^{k+1}-\frac{1}{\rho}(Mx^{k+1}+q^k)\Big), \end{equation}
		 Combining this with the nonexpansiveness of $P_C(.)$ \cite[Corollary~2.4, p.~10]{KS1980} yields
	\begin{equation*}
	\begin{array}{rl}
	& \|x^{k+1}-P_C(x^{k+1}-\frac{1}{\rho}(Qx^{k+1}+q))\|\\
	&\le \|P_C\Big(x^{k+1}-\dfrac{1}{\rho}(Mx^{k+1}+q^k)\Big)-P_C\Big(x^{k+1}-\dfrac{1}{\rho}(Qx^{k+1}+q)\Big)\|\\ & \le \|[x^{k+1}-\dfrac{1}{\rho}(Mx^{k+1}+q^k)]-[x^{k+1}-\dfrac{1}{\rho}(Qx^{k+1}+q)]\|\\
	& = \|[x^{k+1}-\dfrac{1}{\rho}(Qx^{k+1}+\rho x^{k+1}+q-\rho x^k)]- [x^{k+1}-\dfrac{1}{\rho}(Qx^{k+1}+q)]\|\\
	& =\|x^{k+1}-x^k\|<\varepsilon.
	\end{array}
	\end{equation*} Hence \eqref{ineq1n} is valid and, in addition, we have
	$$ \|x^{k+1}-P_C\Big(x^{k+1}-\frac{1}{\rho}(Qx^{k+1}+q)\Big)\|\leq \|x^{k+1}-x^k\|.$$
	From this and \eqref{ineq1} it follows that
	\begin{equation} \label{ineq1nn}
	d(x^{k+1},C^*)\le \ell  \|x^{k+1}-x^k\|\quad\forall k\geq k_0.
	\end{equation} 
		 		Since $C^*$ is closed and nonempty, for each $k\in\{0,1,2,\dots\}$  we can find $y^k\in C^*$ such that
	$d(x^k,C^*)=\|x^k-y^k\|$. Then \eqref{ineq1nn} implies that
	\begin{equation}\label{ineq1_new}
	\|x^{k+1}-y^{k+1}\|\le\ell \|x^{k+1}-x^k\|\quad \forall k\geq k_0.
	\end{equation}
	So, as consequence of \eqref{assym_stab},
	\begin{equation} \label{fact1}
	\lim_{k\to\infty}\|y^{k+1}-x^{k+1}\| =0.
	\end{equation}
Since
	\begin{equation*}
	\|y^{k+1}-y^k\|\le \|y^{k+1}-x^{k+1}\|+\|x^{k+1}-x^{k}\|+ \|x^k-y^k\|,
	\end{equation*}
	it follows that
	\begin{equation} \label{fact2}
	\lim_{k\to\infty}\|y^{k+1}-y^k\| =0.
	\end{equation}
	Let $C_1, C_2, \cdots, C_r$ be the connected components of $C^*$. By
	Lemma \ref{Lemma 3} and \eqref{fact2}, there exist
	$i_0\in\{1,\ldots,r\}$ and $k_1\ge k_0$ such that $y^k\in
	C_{i_0}$ for every $k\ge k_1$. Hence, according to the third assertion of Lemma
	\ref{Lemma 3},
	\begin{equation} \label{fact3}
	f(y^k)=c \quad\forall k\ge k_1
	\end{equation}
	for some  $c\in\mathbb{R}$.
	
	Since \eqref{QP problem} has a solution, by Theorem \ref{basic properties of DCA sequences} we can find a real value $f_*$ such that $\displaystyle \lim_{k\to\infty}f(x^k)=f_*$. 
	 
	By the classical Mean Value Theorem and by the formula $\nabla f(x)=Qx+q$, for
	every $k$ there is $z^k\in (x^k,y^k):=\{(1-t)x^k+ ty^k\, :\,0<t<1\}$
	such that
	$$ f(y^k)-f(x^k)=\langle Qz^k+q, y^k-x^k\rangle.$$
	Since $y^k$ is a KKT point, it holds that
	$0\leq \langle Qy^k+q, x^k-y^k\rangle.$ Adding this inequality and the preceding equality, we get
	\begin{equation}\label{ineq2}
	\begin{array}{rl}
	f(y^k)-f(x^k)&\le\langle Q(z^k-y^k), y^k-x^k\rangle\\
	&\le \|Q\|\|z^k-y^k\||y^k-x^k\|\\
	&\le\|Q\|\|y^k-x^k\|^2.
	\end{array}
	\end{equation}
	
	On one hand, from
\eqref{fact3} and \eqref{ineq2} it follows that
	\begin{equation*}
	c=f(y^k)\le f(x^k)+\|Q\|\,\|y^k-x^k\|^2.
	\end{equation*} As $\displaystyle\lim_{k\to\infty}\left[f(x^k)+\|Q\|\,\|y^k-x^k\|^2\right]=f_*$ due to \eqref{fact1}, this forces
		\begin{equation}\label{ineq3}
		c\leq f_*.
		\end{equation}  
		
On the other hand, since $x^{k+1}=P_C\Big(x^{k+1}-\dfrac{1}{\rho}(Mx^{k+1}+q^k)\Big)$ by \eqref{fixed_point}, the characterization of the
	metric projection on a closed convex set \cite[Theorem~2.3, p.~9]{KS1980} gives us
	$$ \Big\langle\Big[x^{k+1}-\frac{1}{\rho}(Mx^{k+1}+q^k)\Big]-x^{k+1},y- x^{k+1}\Big\rangle\le 0\quad \forall y\in C.$$
Therefore,
	$$ \left\langle Mx^{k+1}+q^k, y^{k+1}-x^{k+1}\right\rangle\ge 0\quad \forall k\in {\mathbb N}.$$
	From this and \eqref{ineq1_new} we get
	\begin{equation*}
	\begin{array}{rl}
	&\langle My^{k+1}+q^k, x^{k+1}-y^{k+1}\rangle\\ & \le \langle My^{k+1}+q^k, x^{k+1}-y^{k+1}\rangle +\langle Mx^{k+1}+q^k, y^{k+1}-x^{k+1}\rangle\\
	&=\langle M(y^{k+1}-x^{k+1}), x^{k+1}-y^{k+1}\rangle\\
	&\le\|M\|\|y^{k+1}-x^{k+1}\|^2\\
	&\le\ell^2\|M\|\|x^{k+1}-x^k\|^2
		\end{array}
	\end{equation*}
	for all $k\geq k_0$. So, setting $\alpha=\ell^2\|M\|$, we have
	\begin{equation}\label{ineq_alpha}
	\langle My^{k+1}+q^k, x^{k+1}-y^{k+1}\rangle\\ 
	\le \alpha\|x^{k+1}-x^k\|^2.
	\end{equation}
		 For each $k\geq k_1$, since $M=Q+\rho I$ and $q^k=q-\rho x^k$, invoking \eqref{ineq_alpha} and using
	\eqref{ineq1_new} once more, we have
	\begin{equation*}\label{ineq4}
	\begin{array}{rl}
	f(x^{k+1})-c&=f(x^{k+1})-f(y^{k+1})\\
	& \le\frac{1}{2}\langle Qx^{k+1},x^{k+1}\rangle +\langle q,x^{k+1}\rangle - \frac{1}{2}\langle Qy^{k+1},y^{k+1}\rangle -\langle q,y^{k+1}\rangle \\
	&=\langle My^{k+1}+q^k, x^{k+1}-y^{k+1}\rangle+\frac{1}{2}\langle Q(x^{k+1}-y^{k+1}),x^{k+1}-y^{k+1}\rangle \\
	&\quad +\rho \langle x^k - y^{k+1},x^{k+1}-y^{k+1}\rangle\\
	
	&=\langle My^{k+1}+q^k, x^{k+1}-y^{k+1}\rangle+\frac{1}{2}\langle Q(x^{k+1}-y^{k+1}),x^{k+1}-y^{k+1}\rangle \\
	&\quad +\rho \langle x^k - x^{k+1},x^{k+1}-y^{k+1}\rangle+\rho \langle x^{k+1}-y^{k+1},x^{k+1}-y^{k+1}\rangle\\
	
	&\le\alpha\|x^{k+1}-x^k\|^2+\frac{1}{2}\|Q\|\|x^{k+1}-y^{k+1}\|^2+\rho \|x^{k+1}-x^k\|\|x^{k+1}-y^{k+1}\|\\
	&\quad +\rho \|x^{k+1}-y^{k+1}\|^2\\
	&\le\left[\alpha + \frac{1}{2}\|Q\|\ell^2+\rho\ell(1+\ell)\right] \|x^{k+1}-x^k\|^2.
	\end{array}
	\end{equation*}
	Therefore, with $\beta:=\alpha + \frac{1}{2}\|Q\|\ell^2+\rho\ell(1+\ell)$, we get
		\begin{equation}\label{ineq4n} f(x^{k+1})\le c+\beta\|x^{k+1}-x^k\|^2.
		\end{equation}
		Letting $k\to\infty$, from \eqref{ineq4n} we can deduce that
	$$f_*=\lim_{k\to\infty} f(x^{k+1})\le c.$$
	
Combining the last expression with \eqref{ineq3} yields $f_*=c$.
	Therefore, by \eqref{ineq4n} and the first assertion of Theorem \ref{basic properties of DCA sequences} we obtain
	$$f(x^{k+1})-f_*\le \beta\|x^{k+1}-x^k\|^2\le \frac{2\beta}{\lambda_1(Q_1)+\lambda_1(Q_2)}(f(x^k)-f(x^{k+1})),$$ where $Q_1=Q+\rho I$ and $Q_2=\rho I$. As $\rho> -\lambda_1(Q)$, putting $\gamma=\lambda_1(Q_1)+\lambda_1(Q_2)$, we see that $\gamma =(\lambda_1(Q)+\rho)+\rho>0$. Therefore, $$f(x^{k+1})-f_*\le \frac{2\beta}{\gamma}\left[\big(f(x^k)-f_*\big)-\big(f(x^{k+1})-f_*\big)\right].$$ Hence
	$$f(x^{k+1})-f_*\le \frac{2\beta}{2\beta+\gamma}(f(x^{k})-f_*).$$
So we have
	$$|f(x^{k+1})-f_*|\le \mu_0 |f(x^{k})-f_*|\quad  \forall\, k\ge k_1,$$
	where $\mu_0:=\frac{2\beta}{2\beta+\gamma}\in (0,1)$. Thus,
	$$|f(x^{k})-f_*|\le \mu_0^{k-k_1}|f(x^{k_1})-f_*|\quad \forall\, k> k_1,$$
	or
	$$|f(x^{k})-f_*|\le r_0\, \mu^{2k}\quad \forall\, k> k_1,$$
	where $r_0:=\mu_0^{-k_1}|f(x^{k_1})-f_*|$ and $\mu:=\mu_0^{1/2}$.
	Hence,
	\begin{equation*}
	\begin{array}{rl}
	|f(x^{k+1})-f(x^k)|&\le |f(x^{k+1})-f_*|+|f(x^{k})-f_*|\\
	&\le r_0\, \mu^{2k+2}+r_0\, \mu^{2k}=r_1\mu^{2k}\quad \forall k>k_1,
	\end{array}
	\end{equation*}
	where $r_1:=r_0(\mu^2+1)$. Consequently, using the first assertion of Theorem \ref{basic properties of DCA sequences} once more, we see that
	$$\|x^{k+1}-x^k\|^2\le \frac{2}{\gamma}(f(x^k)-f(x^{k+1}))\le \frac{2r_1}{\gamma}\mu^{2k}\quad \forall k> k_1.$$
	Thus
	$$\|x^{k+1}-x^k\|\le r\,\mu^{k}\quad \forall k> k_1,$$
	where $r:=\big(\frac{2r_1}{\gamma}\big)^{\frac{1}{2}}$ and $\mu\in (0,1)$. Let $\varepsilon>0$ be given arbitrarily. For each positive integer $p$, we have
	\begin{equation*}
	\begin{array}{rl}
	\|x^{k+p}-x^k\|&\le\|x^{k+p}-x^{k+p-1}\|+\cdots+\|x^{k+1}-x^k\|\\
	&\le r\, \mu^{k+p-1}+\cdots+ r \mu^{k}\\
	&=r\dfrac{1-\mu^p}{1-\mu}\mu^k\le \dfrac{r}{1-\mu}\, \mu^k<\varepsilon,
	\end{array}
	\end{equation*}
	provided that $k$ is large enough. Hence $\{x^k\}$ is a Cauchy
	sequence, and we may assume that it converges to a point $x^*\in C$. By 
	the third assertion of Theorem \ref{basic properties of DCA sequences},
	$x^*\in C^*$. Moreover, passing the inequality
	$$\|x^{k+p}-x^k\|\le \frac{r}{1-\mu}\, \mu^k $$
	to the limit as $p\to\infty$, we get
	$$\|x^k-x^*\|\le \frac{r}{1-\mu}\, \mu^k$$
	for all $k$ large enough. So,
	$$\|x^k-x^*\|^{1/k}\le \left(\frac{r}{1-\mu}\right)^{1/k}\, \mu $$
	for all $k$ large enough. Therefore,
	$$\limsup_{k\to\infty}\|x^k-x^*\|^{1/k}\le\mu<1.$$
	This proves that $\{x^k\}$ converges $R$-linearly to a KKT point of \eqref{QP problem}. $\hfill\Box$
	
\section{Asymptotical Stability of the Algorithm}
We will prove that DCA sequences generated by Algorithm B converge to a locally unique solution of \eqref{QP problem} if the initial points are taken from a suitably-chosen neighborhood of it.

\medskip
First, we have to recall a stability concept that works for discrete dynamical system. Consider \textit{an iteration  algorithm} which generates a unique point $x^{k+1}$, provided that the preceding iteration point $x^k$, $k\in\{0,1,2,\dots\}$, has been defined. Following Leong and Goh \cite[Definition~2]{LeongGoh13}, we can present the concept of \textit{asymptotical stability} of a KKT point as follows.

\begin{Definition}\label{locally unique solution} {\rm The KKT point $x^*$ of \eqref{QP problem}
		is:
		\begin{itemize}
			\item[(i)] \textit{stable} (in the sense of Lyapunov) w.r.t. the iteration algorithm if for any given $\varepsilon > 0$ there exists $\delta δ> 0$ such that whenever
			$x^0\in B(x^*,\delta)$, the DCA sequence generated by the iteration algorithm and the initial point $x^0$ has the property $x^k\in B(x^*,\varepsilon )$ for all $k\geq 0$;
			\item[(ii)] \textit{attractive} if there exists $\delta δ> 0$ such that whenever $x^0\in B(x^*,\delta)$, the DCA sequence generated by the iteration algorithm and the initial point $x^0$ has the property  $\displaystyle\lim_{k\to\infty} x_k = x^*$;
			\item[(iii)] \textit{asymptotically stable} w.r.t. the iteration algorithm  if it is stable and attractive w.r.t. to that algorithm.
		\end{itemize}}
\end{Definition} 

As usual, for an optimization problem $\min\{g(x)\,:\, x\in \Omega\}$ with  $g:\mathbb R^n\rightarrow \mathbb R$ and $\Omega\subset\mathbb R^n$ being respectively a real function and an arbitrary subset, one says that $x^*\in\Omega$ is a {\it locally unique solution} 
		of
		if there exists $\varepsilon>0$ such that $$g(x)>g(x^*)\quad \forall x\in(\Omega\cap B(x^*,\varepsilon))\setminus \{x^*\}.$$
		
		We will need next two lemmas expressing some well-known facts.
		
\begin{Lemma}\label{locally unique solution_strict minimality} {\rm (See, e.g., \cite[Theorem 3.8]{LeeTamYen_book})} {\it If $x^*\in C$ is a locally unique solution of {\rm (\ref{QP problem})},  then there exist $\mu>0$ and $\eta>0$ such that \begin{eqnarray}\label{property (1)}f(x)-f(x^*)\geq \eta \|x-x^*\|^2\quad \mbox{for\ every}\ \, x\in C\cap B(x^*,\mu).\end{eqnarray}}
\end{Lemma}

\begin{Lemma}\label{Objective function on a KKT interval} {\rm (See, e.g.,  \cite[Proof of Lemma 4]{Bomze_Danninger_1994} and \cite[Lemma~1]{ATY2})} If the KKT point set $ C^*$ contains a segment $[u,x]$, then the restriction of $f$ on that segment is a constant function.
\end{Lemma}

The main result of this section can be formulated as follows.

\begin{Theorem}\label{Convergence near local unique solutions} Consider Algorithm B and require additionally that $\rho >\|Q\|$. Suppose $x^*$ is  a  locally unique solution of problem {\rm (\ref{QP problem})}. In that case, for any  $\varepsilon>0$ there exists $\delta>0$ such that if $x^0\in C\cap B(x^*,\delta)$ and if $\{x^k\}$ is the DCA sequence generated by Algorithm B and the initial point $x^0$, then
	\begin{itemize}
		\item[{\rm (a)}] $x^k\in C\cap B(x^*,\varepsilon)\,$ for any $k\geq 0$; 
		\item[{\rm (b)}] $x^k\to x^*\,$ as $k\to\infty$. 
		\end{itemize}
		In other words,  $x^*$ is asymptotically stable w.r.t. Algorithm B.
\end{Theorem}
\noindent \textit{Proof.}  Suppose that $\rho >\|Q\|$ and $x^*$  is a locally unique solution of (\ref{QP problem}). By Lemma~\ref{locally unique solution_strict minimality} we can select constants $\mu>0$ and $\eta>0$ such that \eqref{property (1)} holds. For any given $\varepsilon>0$, by replacing $\varepsilon$ with a smaller one (if necessary), we may assume that $\varepsilon\in (0,\mu)$ and $\varepsilon<\mu(1-\rho^{-1}\|Q\|).$ Since 
$$f(x)-f(x^*)>0\quad \forall x\in \big(C\cap B(x^*,\varepsilon)\big)\setminus\{x^*\}$$ by \eqref{property (1)}, the continuity of $f$ implies the existence of $\delta\in (0,\mu)$ satisfying 
\begin{eqnarray}\label{property (2)}\displaystyle\frac{1}{\eta^{1/2}}\big(f(x)-f(x^*)\big)^{1/2}<\varepsilon \quad \forall x\in C\cap B(x^*,\delta).\end{eqnarray} 

First, let us show that the assertion about stability in the sense of Lyapunov of DCA sequences generated by Algorithm B is valid for the chosen number $\delta>0$. Fix any $x^0\in C\cap B(x^*,\delta)$. As $\delta<\varepsilon$, for $k=0$ we have $x^k\in C\cap B(x^*,\varepsilon)$. To proceed by induction, suppose that the last inclusion holds for some $k\geq 0$. Since $x^*$ is a locally unique solution of  (\ref{QP problem}), it is a KKT point of that problem, i.e., 
\begin{eqnarray}\label{property (3)}\big(Qx^*+q\big)^T(x-x^*)\geq 0\quad \forall x\in C.\end{eqnarray} It follows that
\begin{eqnarray}\label{property (4)}x^*=P_C\big(x^*-\frac{1}{\rho}(Qx^*+q)\big).\end{eqnarray} Indeed, by the characterization of the metric projection \cite[Theorem~2.3, p.~9]{KS1980}, (\ref{property (4)}) is valid if and only if 
$$\Big(\big[x^*-\frac{1}{\rho}(Qx^*+q)\big]-x^*\Big)^T(x-x^*)\leq 0\quad \forall x\in C.$$
The latter is equivalent to (\ref{property (3)}). Using \eqref{Iteration formula of Algorithm B}, (\ref{property (4)}), and the nonexpansiveness of the metric projection \cite[Corollary~2.4, p.~10]{KS1980}, we have 
\begin{eqnarray*} \|x^{k+1}-x^*\|&=& \big\|P_C\big(x^{k+1}-\frac{1}{\rho}(Mx^{k+1}+q^k)\big)-P_C\big(x^*-\frac{1}{\rho}(Qx^*+q)\big)\big\|\\
	&\leq& \big\|\big[x^{k+1}-\frac{1}{\rho}(Mx^{k+1}+q^k)\big]-\big[x^*-\frac{1}{\rho}(Qx^*+q)\big]\big\|\\
	& =& \big\|\big[x^{k+1}-\frac{1}{\rho}\big((\rho I + Q)x^{k+1}+q-\rho x^k\big)\big]-\big[x^*-\frac{1}{\rho}(Qx^*+q)\big]\big\|\\
	& =& \big\|(x^k-x^*)+\frac{1}{\rho}Q(x^*-x^{k+1})\big\|\\
	&\leq& \|x^k-x^*\|+\frac{1}{\rho}\|Q\|\|x^*-x^{k+1}\|.\end{eqnarray*}
	Then we obtain
	\begin{eqnarray*}
	\|x^{k+1}-x^*\| \leq (1-\frac{1}{\rho}\|Q\|)^{-1}\|x^k-x^*\|
	\leq (1-\frac{1}{\rho}\|Q\|)^{-1}\varepsilon<\mu,
\end{eqnarray*} where the strict inequality follows from the property  $\varepsilon<\mu(1-\rho^{-1}\|Q\|).$ Thus, $x^{k+1}\in C\cap B(x^*,\mu)$. Applying (\ref{property (1)}) and the inequality $f(x^k)\geq f(x^{k+1})$ which holds for any $k\geq 0$ (see Remark \ref{Remark 1}), we get 
\begin{eqnarray*} \|x^{k+1}-x^*\|^2&\leq & \frac{1}{\eta}\big(f(x^{k+1})-f(x^*)\big)\\
	&\leq & \frac{1}{\eta}\big(f(x^k)-f(x^*)\big)\\
	& & \vdots\\
	&\leq & \frac{1}{\eta}\big(f(x^0)-f(x^*)\big).
\end{eqnarray*} Hence,
$$\|x^{k+1}-x^*\|\leq\displaystyle\frac{1}{\eta^{1/2}}\big(f(x^0)-f(x^*)\big)^{1/2}.$$ Since $x^0\in C\cap B(x^*,\delta)$, combining this with (\ref{property (2)}) we obtain $\|x^{k+1}-x^*\|<\varepsilon$ which means that $x^{k+1}\in C\cap B(x^*,\varepsilon)$. Thus, we have proved that $x^k\in C\cap B(x^*,\varepsilon)$ for every $k\geq 0$.

Next, to obtain the assertion about the attractiveness of DCA sequences generated by Algorithm B, we observe by the just obtained stability result that for any $\varepsilon>0$ there exists $\delta=\delta(\varepsilon)>0$ such that if $x^0\in C\cap B(x^*,\delta)$ and if $\{x^k\}$ is the DCA sequence generated by Algorithm B and the initial point $x^0$, then the property in (a) is valid. Without loss of generality, we may assume that $\varepsilon\in (0,\mu)$ and $\delta\in (0,\varepsilon)$. By taking a smaller positive $\varepsilon>0$ and choosing the corresponding $\delta=\delta(\varepsilon)$ such that the property in (a) is valid, we can have the following:  If $x^0\in C\cap B(x^*,\delta)$ and if $\{x^k\}$ is the DCA sequence generated by Algorithm B and the initial point $x^0$, then the property in (b) holds. Indeed, if this claim was false, we would find sequences $\varepsilon_j\to 0^+$ and $\delta_j\to 0^+$ such that for each $j\in\mathbb N$ we have $\varepsilon_j\in (0,\mu)$,  $\delta_j\in (0,\varepsilon_j)$,  and the stability assertion is valid for the pair $(\delta,\varepsilon):=(\delta_j,\varepsilon_j)$. Moreover, for each $j$, there exists some  $x^{0,j}\in C\cap B(x^*,\delta_j)$ such that the DCA sequence  $\{x^{k,j}\}$ generated by Algorithm B and the initial point $x^{0,j}$ does not converge to $x^*$. Then we can select a subsequence of $\{x^{k,j}\}$ which converges to a point  \begin{eqnarray}\label{property (5a)}\widetilde x^j\in C\cap \overline{B}(x^*,\varepsilon_j)\subset C\cap B(x^*,\mu),\end{eqnarray} where $\widetilde x^j\neq x^*$. By Theorem \ref{basic properties of DCA sequences} we have $\widetilde x^j\in C^*$ for $j=1,2,\dots.$ Observe that
\begin{eqnarray}\label{property (5)}\displaystyle\lim_{j\to\infty}\widetilde x^j=x^*.\end{eqnarray} 
For each $j$, one can find a natural number $k(j)\geq 1$ such that $\varepsilon_{j+k(j)}<\|\widetilde x^j-x^*\|$.  Then, by \eqref{property (5a)} one has $$\|\widetilde x^{j+k(j)}-x^*\|<\|\widetilde x^j-x^*\|.$$ Choose $z^1=\widetilde x^1$ and set $z^{p+1}:=\widetilde x^{p+k(p)}$ for $p=1,2,\dots$. It is clear that $\{z^p\}$ is a subsequence of $\{\widetilde x^j\}$ and  $z^p\neq z^{p'}$ whenever $p'\neq p$. Hence, by considering a subsequence (if necessary), we can assume that $\widetilde x^j\neq \widetilde x^\ell$ whenever $j\neq\ell$. Since the number of pseudo-faces of $C$ is finite, by~\eqref{property (5)} there must exists an index set $\alpha\subset\{1,\dots,m\}$ such that the pseudo-face
$$F_\alpha:=\{x\in\mathbb R^n\,:\,  A_\alpha x=b_\alpha,\ A_{\bar \alpha}x>b_{\bar \alpha}\}$$ of $C$
contains infinite number of the members of the sequence $\{\widetilde x^j\}$. Without loss of generality, we may assume that the whole sequence $\{\widetilde x^j\}$ is contained in $F_\alpha$. By \cite[Lemma 4.1]{LeeTamYen_book},  the intersection $ C^*\cap F_\alpha$ is a convex set. Hence, according to Lemma~\ref{Objective function on a KKT interval}, the restriction of $f$ on  $C^*\cap F_\alpha$ is a constant function. Using (\ref{property (5)}), from this we can deduce that the equality $f(\widetilde x^j)=f(x^*)$ holds for all $j$. As $\widetilde x^j\neq x^*$ for every $j$, the last equality contradicts (\ref{property (1)}). Our claim has been proved. $\hfill\square$  

\section{Further Analysis}

In this final section, we will analyze the influence of the decomposition parameter $\rho$ for the rates of convergence of the algorithms A and B. We also compare the effectiveness of Algorithm~B with that of Algorithm~A. These algorithms were implemented in the Visual C++ 2010 environment, and performed on a PC Intel Core${^\textit{TM}}$ i7 (4 x 2.0 GHz) processor, 4GB RAM. The CPLEX 11.2 solver is used to solve linear and convex quadratic problems.

\medskip 
Recall that, for Algorithm A, the parameter $\rho>0$ has to satisfy the inequality $\rho\geq\lambda_n(Q)$. For Algorithm B,  $\rho>0$ must satisfy the strict inequality $\rho>-\lambda_1(Q)$. 

\medskip
We now present the results of our tests in using the algorithms A and B to solve problem \eqref{QP problem} for the dimensions $n=10$, $n=20$, $n=40$, $n=60$, $n=80$. With $\beta_i\in [0,10]$ for $i=1,\dots,n$ being generated randomly, the following two types of constraint sets have been considered:
$$C=\Big\{x\in\mathbb R^n\; :\; x\geq 0,\ ix_i\geq \beta_i,\ i=1,\dots, n,\  \sum_{i=1}^n ix_i\leq 5000\Big\}$$
and 
$$C=\Big\{x\in\mathbb R^n\; :\; x\geq 0,\ ix_i\geq \beta_i,\ i=1,\dots, n,\  10\leq x_1+\sum_{i=2}^n0.1 ix_i\leq 100\Big\}.$$ It is easy to express each of these sets as the solution set of the linear inequality system $Ax\geq b$ with a suitably chosen matrix  $A\in\mathbb  R^{m\times n}$ and a vector $b\in\mathbb  R^m$. 
Fixing a dimension $n\in\{10, 20, 40, 60, 80\}$, we generate randomly a symmetric matrix $Q\in\mathbb R^{n\times n}$ and a vector $q\in\mathbb  R^n$ with the requirement that all their components belong to the segment $[0,10]$. The initial point $x^0\in \mathbb R^{n\times n}$ is generated randomly with the requirement that all its components belong to the segment $[0,5]$. Then, we start testing Algorithm A with  $\rho=\lambda_n(Q)$ if $\lambda_n(Q)>0$ and $\rho=0.1$ otherwise. For our convenience, this $\rho$ is called \textit{the smallest decomposition parameter} for  Algorithm A. Similarly, we start testing Algorithm B with  $\rho=-\lambda_1(Q)+0.1$ if $\lambda_1(Q)<0$ and $\rho=0.1$ otherwise. This $\rho$ is said to be \textit{the smallest decomposition parameter} for  Algorithm~B. The stopping criterion is $\|x^{k+1}-x^k\|\leq 10^{-6}$ and the allowed largest number of steps is 1000. After testing Algorithm A (resp., Algorithm B) for a decomposition parameter $\rho$, we increase $\rho$ by 1.5 times and let the algorithm to run again. 

\medskip
In Table 1, the second rows of the tables a) and b) correspond to the smallest decomposition parameters for Algorithm A and Algorithm B, respectively. The decomposition parameters of the test reported in the third rows are 1.5 times of the smallest decomposition parameters. The decomposition parameters of the test reported in the fourth rows are 1.5 times of the just mentioned decomposition parameters; and so on... In the tables a) and b), the first column presents the ordinal number of the tests.  The second one indicates  the numbers of iterations. The third one reports  the running times.  And the fourth column contains  the decomposition parameters. There are only 11 records in table a) because for larger decomposition parameters, the numbers of steps exceed 1000. For the same reason, table b) just contains 18 records. 

\medskip
The contents of Tables 2--6 are similar to those of Table 1.

\medskip
 With any  $n$ belonging to the set $\{10, 20, 40, 60, 80\}$, a careful analysis of these Tables allows us to observe that:

\textit{$\bullet$ For both algorithms, if $\rho$ increases, then the running time, as well as the number of computation steps, increases;} 

\textit{$\bullet$ For the rows of the tables a) and b) with the same ordinal number, Algorithm B is much more efficient than Algorithm A (for example, the running time of the first one is much smaller than that of the second one).}

\medskip
Due to the space limitation, we only present the test results for $n=10,\, 40,\, 80$.

% Table generated by Excel2LaTeX from sheet 'FileTestTG10'
\begin{table}[htbp]
	\centering
	\caption{The test results for \textit{n} = 10 with the 1st type constraint}
	\begin{tabular}{rrrrrrrrcr}\\
		\cline{1-4}\cline{7-10}    \multicolumn{1}{|c|}{No.} & \multicolumn{1}{c|}{Step} & \multicolumn{1}{c|}{Time} & \multicolumn{1}{c|}{roA} &       & \multicolumn{1}{c|}{} & \multicolumn{1}{c|}{No.} & \multicolumn{1}{c|}{Step} & \multicolumn{1}{c|}{Time} & \multicolumn{1}{c|}{roB} \\
		\cline{1-4}\cline{7-10}    \multicolumn{1}{|r|}{1} & \multicolumn{1}{r|}{5} & \multicolumn{1}{r|}{0.239} & \multicolumn{1}{r|}{48.802} &       & \multicolumn{1}{r|}{} & \multicolumn{1}{r|}{1} & \multicolumn{1}{r|}{4} & \multicolumn{1}{r|}{0.127} & \multicolumn{1}{r|}{9.380} \\
		\cline{1-4}\cline{7-10}    \multicolumn{1}{|r|}{2} & \multicolumn{1}{r|}{12} & \multicolumn{1}{r|}{0.222} & \multicolumn{1}{r|}{73.203} &       & \multicolumn{1}{r|}{} & \multicolumn{1}{r|}{2} & \multicolumn{1}{r|}{4} & \multicolumn{1}{r|}{0.125} & \multicolumn{1}{r|}{14.070} \\
		\cline{1-4}\cline{7-10}    \multicolumn{1}{|r|}{3} & \multicolumn{1}{r|}{22} & \multicolumn{1}{r|}{0.274} & \multicolumn{1}{r|}{109.805} &       & \multicolumn{1}{r|}{} & \multicolumn{1}{r|}{3} & \multicolumn{1}{r|}{5} & \multicolumn{1}{r|}{0.114} & \multicolumn{1}{r|}{21.105} \\
		\cline{1-4}\cline{7-10}    \multicolumn{1}{|r|}{4} & \multicolumn{1}{r|}{37} & \multicolumn{1}{r|}{0.416} & \multicolumn{1}{r|}{164.707} &       & \multicolumn{1}{r|}{} & \multicolumn{1}{r|}{4} & \multicolumn{1}{r|}{6} & \multicolumn{1}{r|}{0.135} & \multicolumn{1}{r|}{31.658} \\
		\cline{1-4}\cline{7-10}    \multicolumn{1}{|r|}{5} & \multicolumn{1}{r|}{59} & \multicolumn{1}{r|}{0.718} & \multicolumn{1}{r|}{247.060} &       & \multicolumn{1}{r|}{} & \multicolumn{1}{r|}{5} & \multicolumn{1}{r|}{8} & \multicolumn{1}{r|}{0.210} & \multicolumn{1}{r|}{47.487} \\
		\cline{1-4}\cline{7-10}    \multicolumn{1}{|r|}{6} & \multicolumn{1}{r|}{91} & \multicolumn{1}{r|}{0.947} & \multicolumn{1}{r|}{370.590} &       & \multicolumn{1}{r|}{} & \multicolumn{1}{r|}{6} & \multicolumn{1}{r|}{10} & \multicolumn{1}{r|}{0.227} & \multicolumn{1}{r|}{71.231} \\
		\cline{1-4}\cline{7-10}    \multicolumn{1}{|r|}{7} & \multicolumn{1}{r|}{139} & \multicolumn{1}{r|}{1.364} & \multicolumn{1}{r|}{555.886} &       & \multicolumn{1}{r|}{} & \multicolumn{1}{r|}{7} & \multicolumn{1}{r|}{13} & \multicolumn{1}{r|}{0.296} & \multicolumn{1}{r|}{106.846} \\
		\cline{1-4}\cline{7-10}    \multicolumn{1}{|r|}{8} & \multicolumn{1}{r|}{210} & \multicolumn{1}{r|}{2.050} & \multicolumn{1}{r|}{833.829} &       & \multicolumn{1}{r|}{} & \multicolumn{1}{r|}{8} & \multicolumn{1}{r|}{17} & \multicolumn{1}{r|}{0.419} & \multicolumn{1}{r|}{160.269} \\
		\cline{1-4}\cline{7-10}    \multicolumn{1}{|r|}{9} & \multicolumn{1}{r|}{316} & \multicolumn{1}{r|}{3.019} & \multicolumn{1}{r|}{1250.743} &       & \multicolumn{1}{r|}{} & \multicolumn{1}{r|}{9} & \multicolumn{1}{r|}{24} & \multicolumn{1}{r|}{0.576} & \multicolumn{1}{r|}{240.404} \\
		\cline{1-4}\cline{7-10}    \multicolumn{1}{|r|}{10} & \multicolumn{1}{r|}{474} & \multicolumn{1}{r|}{4.593} & \multicolumn{1}{r|}{1876.114} &       & \multicolumn{1}{r|}{} & \multicolumn{1}{r|}{10} & \multicolumn{1}{r|}{34} & \multicolumn{1}{r|}{0.787} & \multicolumn{1}{r|}{360.606} \\
		\cline{1-4}\cline{7-10}    \multicolumn{1}{|r|}{11} & \multicolumn{1}{r|}{710} & \multicolumn{1}{r|}{7.006} & \multicolumn{1}{r|}{2814.171} &       & \multicolumn{1}{r|}{} & \multicolumn{1}{r|}{11} & \multicolumn{1}{r|}{49} & \multicolumn{1}{r|}{1.312} & \multicolumn{1}{r|}{540.909} \\
		\cline{1-4}\cline{7-10}          &       &       &       &       & \multicolumn{1}{r|}{} & \multicolumn{1}{r|}{12} & \multicolumn{1}{r|}{72} & \multicolumn{1}{r|}{1.775} & \multicolumn{1}{r|}{811.363} \\
		\cline{7-10}          &       & \multicolumn{1}{c}{a)} &       &       & \multicolumn{1}{r|}{} & \multicolumn{1}{r|}{13} & \multicolumn{1}{r|}{106} & \multicolumn{1}{r|}{2.921} & \multicolumn{1}{r|}{1217.044} \\
		\cline{7-10}          &       &       &       &       & \multicolumn{1}{r|}{} & \multicolumn{1}{r|}{14} & \multicolumn{1}{r|}{157} & \multicolumn{1}{r|}{4.244} & \multicolumn{1}{r|}{1825.567} \\
		\cline{7-10}          &       &       &       &       & \multicolumn{1}{r|}{} & \multicolumn{1}{r|}{15} & \multicolumn{1}{r|}{233} & \multicolumn{1}{r|}{6.155} & \multicolumn{1}{r|}{2738.350} \\
		\cline{7-10}          &       &       &       &       & \multicolumn{1}{r|}{} & \multicolumn{1}{r|}{16} & \multicolumn{1}{r|}{348} & \multicolumn{1}{r|}{9.053} & \multicolumn{1}{r|}{4107.525} \\
		\cline{7-10}          &       &       &       &       & \multicolumn{1}{r|}{} & \multicolumn{1}{r|}{17} & \multicolumn{1}{r|}{520} & \multicolumn{1}{r|}{13.852} & \multicolumn{1}{r|}{6161.288} \\
		\cline{7-10}          &       &       &       &       & \multicolumn{1}{r|}{} & \multicolumn{1}{r|}{18} & \multicolumn{1}{r|}{778} & \multicolumn{1}{r|}{20.276} & \multicolumn{1}{r|}{9241.932} \\
		\cline{7-10}          &       &       &       &       &       &       &       &       &  \\
		&       &       &       &       &       &       &       & b)    &  \\
	\end{tabular}%
	\label{tab:addlabel}%
\end{table}%

%\input{FileTest_10.tex} 

% Table generated by Excel2LaTeX from sheet 'FileTestHT10'
\begin{table}[htbp]
	\centering
	\caption{The test results for \textit{n} = 10 with the 2nd type constraint}
	\begin{tabular}{rrrrrrrrcr}\\
		\cline{1-4}\cline{7-10}    \multicolumn{1}{|c|}{No.} & \multicolumn{1}{c|}{Step} & \multicolumn{1}{c|}{Time} & \multicolumn{1}{c|}{roA} &       & \multicolumn{1}{c|}{} & \multicolumn{1}{c|}{No.} & \multicolumn{1}{c|}{Step} & \multicolumn{1}{c|}{Time} & \multicolumn{1}{c|}{roB} \\
		\cline{1-4}\cline{7-10}    \multicolumn{1}{|r|}{1} & \multicolumn{1}{r|}{3} & \multicolumn{1}{r|}{0.189} & \multicolumn{1}{r|}{47.763} &       & \multicolumn{1}{r|}{} & \multicolumn{1}{r|}{1} & \multicolumn{1}{r|}{3} & \multicolumn{1}{r|}{0.131} & \multicolumn{1}{r|}{15.645} \\
		\cline{1-4}\cline{7-10}    \multicolumn{1}{|r|}{2} & \multicolumn{1}{r|}{7} & \multicolumn{1}{r|}{0.210} & \multicolumn{1}{r|}{71.644} &       & \multicolumn{1}{r|}{} & \multicolumn{1}{r|}{2} & \multicolumn{1}{r|}{4} & \multicolumn{1}{r|}{0.175} & \multicolumn{1}{r|}{23.468} \\
		\cline{1-4}\cline{7-10}    \multicolumn{1}{|r|}{3} & \multicolumn{1}{r|}{13} & \multicolumn{1}{r|}{0.285} & \multicolumn{1}{r|}{107.467} &       & \multicolumn{1}{r|}{} & \multicolumn{1}{r|}{3} & \multicolumn{1}{r|}{4} & \multicolumn{1}{r|}{0.167} & \multicolumn{1}{r|}{35.201} \\
		\cline{1-4}\cline{7-10}    \multicolumn{1}{|r|}{4} & \multicolumn{1}{r|}{21} & \multicolumn{1}{r|}{0.233} & \multicolumn{1}{r|}{161.200} &       & \multicolumn{1}{r|}{} & \multicolumn{1}{r|}{4} & \multicolumn{1}{r|}{6} & \multicolumn{1}{r|}{0.252} & \multicolumn{1}{r|}{52.802} \\
		\cline{1-4}\cline{7-10}    \multicolumn{1}{|r|}{5} & \multicolumn{1}{r|}{33} & \multicolumn{1}{r|}{0.335} & \multicolumn{1}{r|}{241.800} &       & \multicolumn{1}{r|}{} & \multicolumn{1}{r|}{5} & \multicolumn{1}{r|}{7} & \multicolumn{1}{r|}{0.206} & \multicolumn{1}{r|}{79.203} \\
		\cline{1-4}\cline{7-10}    \multicolumn{1}{|r|}{6} & \multicolumn{1}{r|}{51} & \multicolumn{1}{r|}{0.527} & \multicolumn{1}{r|}{362.700} &       & \multicolumn{1}{r|}{} & \multicolumn{1}{r|}{6} & \multicolumn{1}{r|}{9} & \multicolumn{1}{r|}{0.329} & \multicolumn{1}{r|}{118.805} \\
		\cline{1-4}\cline{7-10}    \multicolumn{1}{|r|}{7} & \multicolumn{1}{r|}{77} & \multicolumn{1}{r|}{0.729} & \multicolumn{1}{r|}{544.049} &       & \multicolumn{1}{r|}{} & \multicolumn{1}{r|}{7} & \multicolumn{1}{r|}{12} & \multicolumn{1}{r|}{0.298} & \multicolumn{1}{r|}{178.207} \\
		\cline{1-4}\cline{7-10}    \multicolumn{1}{|r|}{8} & \multicolumn{1}{r|}{115} & \multicolumn{1}{r|}{1.029} & \multicolumn{1}{r|}{816.074} &       & \multicolumn{1}{r|}{} & \multicolumn{1}{r|}{8} & \multicolumn{1}{r|}{16} & \multicolumn{1}{r|}{0.506} & \multicolumn{1}{r|}{267.310} \\
		\cline{1-4}\cline{7-10}    \multicolumn{1}{|r|}{9} & \multicolumn{1}{r|}{171} & \multicolumn{1}{r|}{1.802} & \multicolumn{1}{r|}{1224.111} &       & \multicolumn{1}{r|}{} & \multicolumn{1}{r|}{9} & \multicolumn{1}{r|}{22} & \multicolumn{1}{r|}{0.830} & \multicolumn{1}{r|}{400.966} \\
		\cline{1-4}\cline{7-10}    \multicolumn{1}{|r|}{10} & \multicolumn{1}{r|}{255} & \multicolumn{1}{r|}{2.363} & \multicolumn{1}{r|}{1836.167} &       & \multicolumn{1}{r|}{} & \multicolumn{1}{r|}{10} & \multicolumn{1}{r|}{31} & \multicolumn{1}{r|}{1.073} & \multicolumn{1}{r|}{601.449} \\
		\cline{1-4}\cline{7-10}    \multicolumn{1}{|r|}{11} & \multicolumn{1}{r|}{380} & \multicolumn{1}{r|}{3.637} & \multicolumn{1}{r|}{2754.250} &       & \multicolumn{1}{r|}{} & \multicolumn{1}{r|}{11} & \multicolumn{1}{r|}{44} & \multicolumn{1}{r|}{1.043} & \multicolumn{1}{r|}{902.173} \\
		\cline{1-4}\cline{7-10}    \multicolumn{1}{|r|}{12} & \multicolumn{1}{r|}{567} & \multicolumn{1}{r|}{5.133} & \multicolumn{1}{r|}{4131.375} &       & \multicolumn{1}{r|}{} & \multicolumn{1}{r|}{12} & \multicolumn{1}{r|}{65} & \multicolumn{1}{r|}{1.543} & \multicolumn{1}{r|}{1353.259} \\
		\cline{1-4}\cline{7-10}    \multicolumn{1}{|r|}{13} & \multicolumn{1}{r|}{847} & \multicolumn{1}{r|}{7.546} & \multicolumn{1}{r|}{6197.063} &       & \multicolumn{1}{r|}{} & \multicolumn{1}{r|}{13} & \multicolumn{1}{r|}{95} & \multicolumn{1}{r|}{2.628} & \multicolumn{1}{r|}{2029.889} \\
		\cline{1-4}\cline{7-10}          &       &       &       &       & \multicolumn{1}{r|}{} & \multicolumn{1}{r|}{14} & \multicolumn{1}{r|}{141} & \multicolumn{1}{r|}{3.178} & \multicolumn{1}{r|}{3044.833} \\
		\cline{7-10}          &       & \multicolumn{1}{c}{a)} &       &       & \multicolumn{1}{r|}{} & \multicolumn{1}{r|}{15} & \multicolumn{1}{r|}{210} & \multicolumn{1}{r|}{4.914} & \multicolumn{1}{r|}{4567.250} \\
		\cline{7-10}          &       &       &       &       & \multicolumn{1}{r|}{} & \multicolumn{1}{r|}{16} & \multicolumn{1}{r|}{313} & \multicolumn{1}{r|}{7.410} & \multicolumn{1}{r|}{6850.875} \\
		\cline{7-10}          &       &       &       &       & \multicolumn{1}{r|}{} & \multicolumn{1}{r|}{17} & \multicolumn{1}{r|}{467} & \multicolumn{1}{r|}{11.348} & \multicolumn{1}{r|}{10276.313} \\
		\cline{7-10}          &       &       &       &       & \multicolumn{1}{r|}{} & \multicolumn{1}{r|}{18} & \multicolumn{1}{r|}{699} & \multicolumn{1}{r|}{17.063} & \multicolumn{1}{r|}{15414.469} \\
		\cline{7-10}          &       &       &       &       &       &       &       &       &  \\
		&       &       &       &       &       &       &       & b)    &  \\
	\end{tabular}%
	\label{tab:addlabel}%
\end{table}%

% Table generated by Excel2LaTeX from sheet 'FileTestTG40'
\begin{table}[htbp]
	\centering
	\caption{The test results for \textit{n} = 40  with the 1st type constraint}
	\begin{tabular}{rrrrrrrrcr}\\
		\cline{1-4}\cline{7-10}    \multicolumn{1}{|c|}{No.} & \multicolumn{1}{c|}{Step} & \multicolumn{1}{c|}{Time} & \multicolumn{1}{c|}{roA} &       & \multicolumn{1}{c|}{} & \multicolumn{1}{c|}{No.} & \multicolumn{1}{c|}{Step} & \multicolumn{1}{c|}{Time} & \multicolumn{1}{c|}{roB} \\
		\cline{1-4}\cline{7-10}    \multicolumn{1}{|r|}{1} & \multicolumn{1}{r|}{8} & \multicolumn{1}{r|}{0.621} & \multicolumn{1}{r|}{194.883} &       & \multicolumn{1}{r|}{} & \multicolumn{1}{r|}{1} & \multicolumn{1}{r|}{5} & \multicolumn{1}{r|}{0.320} & \multicolumn{1}{r|}{32.917} \\
		\cline{1-4}\cline{7-10}    \multicolumn{1}{|r|}{2} & \multicolumn{1}{r|}{20} & \multicolumn{1}{r|}{0.664} & \multicolumn{1}{r|}{292.324} &       & \multicolumn{1}{r|}{} & \multicolumn{1}{r|}{2} & \multicolumn{1}{r|}{6} & \multicolumn{1}{r|}{0.386} & \multicolumn{1}{r|}{49.375} \\
		\cline{1-4}\cline{7-10}    \multicolumn{1}{|r|}{3} & \multicolumn{1}{r|}{65} & \multicolumn{1}{r|}{1.498} & \multicolumn{1}{r|}{657.729} &       & \multicolumn{1}{r|}{} & \multicolumn{1}{r|}{3} & \multicolumn{1}{r|}{7} & \multicolumn{1}{r|}{0.454} & \multicolumn{1}{r|}{74.062} \\
		\cline{1-4}\cline{7-10}    \multicolumn{1}{|r|}{4} & \multicolumn{1}{r|}{106} & \multicolumn{1}{r|}{2.256} & \multicolumn{1}{r|}{986.594} &       & \multicolumn{1}{r|}{} & \multicolumn{1}{r|}{4} & \multicolumn{1}{r|}{8} & \multicolumn{1}{r|}{0.509} & \multicolumn{1}{r|}{111.094} \\
		\cline{1-4}\cline{7-10}    \multicolumn{1}{|r|}{5} & \multicolumn{1}{r|}{167} & \multicolumn{1}{r|}{3.255} & \multicolumn{1}{r|}{1479.891} &       & \multicolumn{1}{r|}{} & \multicolumn{1}{r|}{5} & \multicolumn{1}{r|}{11} & \multicolumn{1}{r|}{0.670} & \multicolumn{1}{r|}{166.641} \\
		\cline{1-4}\cline{7-10}    \multicolumn{1}{|r|}{6} & \multicolumn{1}{r|}{259} & \multicolumn{1}{r|}{4.925} & \multicolumn{1}{r|}{2219.837} &       & \multicolumn{1}{r|}{} & \multicolumn{1}{r|}{6} & \multicolumn{1}{r|}{15} & \multicolumn{1}{r|}{0.947} & \multicolumn{1}{r|}{249.961} \\
		\cline{1-4}\cline{7-10}    \multicolumn{1}{|r|}{7} & \multicolumn{1}{r|}{397} & \multicolumn{1}{r|}{7.451} & \multicolumn{1}{r|}{3329.755} &       & \multicolumn{1}{r|}{} & \multicolumn{1}{r|}{7} & \multicolumn{1}{r|}{20} & \multicolumn{1}{r|}{1.238} & \multicolumn{1}{r|}{374.941} \\
		\cline{1-4}\cline{7-10}    \multicolumn{1}{|r|}{8} & \multicolumn{1}{r|}{604} & \multicolumn{1}{r|}{11.236} & \multicolumn{1}{r|}{4994.632} &       & \multicolumn{1}{r|}{} & \multicolumn{1}{r|}{8} & \multicolumn{1}{r|}{28} & \multicolumn{1}{r|}{1.734} & \multicolumn{1}{r|}{562.412} \\
		\cline{1-4}\cline{7-10}    \multicolumn{1}{|r|}{9} & \multicolumn{1}{r|}{915} & \multicolumn{1}{r|}{17.078} & \multicolumn{1}{r|}{7491.948} &       & \multicolumn{1}{r|}{} & \multicolumn{1}{r|}{9} & \multicolumn{1}{r|}{40} & \multicolumn{1}{r|}{2.477} & \multicolumn{1}{r|}{843.618} \\
		\cline{1-4}\cline{7-10}          &       &       &       &       & \multicolumn{1}{r|}{} & \multicolumn{1}{r|}{10} & \multicolumn{1}{r|}{57} & \multicolumn{1}{r|}{3.507} & \multicolumn{1}{r|}{1265.427} \\
		\cline{7-10}          &       & \multicolumn{1}{c}{a)} &       &       & \multicolumn{1}{r|}{} & \multicolumn{1}{r|}{11} & \multicolumn{1}{r|}{84} & \multicolumn{1}{r|}{5.061} & \multicolumn{1}{r|}{1898.141} \\
		\cline{7-10}          &       &       &       &       & \multicolumn{1}{r|}{} & \multicolumn{1}{r|}{12} & \multicolumn{1}{r|}{123} & \multicolumn{1}{r|}{7.938} & \multicolumn{1}{r|}{2847.211} \\
		\cline{7-10}          &       &       &       &       & \multicolumn{1}{r|}{} & \multicolumn{1}{r|}{13} & \multicolumn{1}{r|}{182} & \multicolumn{1}{r|}{11.181} & \multicolumn{1}{r|}{4270.817} \\
		\cline{7-10}          &       &       &       &       & \multicolumn{1}{r|}{} & \multicolumn{1}{r|}{14} & \multicolumn{1}{r|}{271} & \multicolumn{1}{r|}{16.625} & \multicolumn{1}{r|}{6406.225} \\
		\cline{7-10}          &       &       &       &       & \multicolumn{1}{r|}{} & \multicolumn{1}{r|}{15} & \multicolumn{1}{r|}{403} & \multicolumn{1}{r|}{24.807} & \multicolumn{1}{r|}{9609.338} \\
		\cline{7-10}          &       &       &       &       & \multicolumn{1}{r|}{} & \multicolumn{1}{r|}{16} & \multicolumn{1}{r|}{602} & \multicolumn{1}{r|}{37.672} & \multicolumn{1}{r|}{14414.006} \\
		\cline{7-10}          &       &       &       &       & \multicolumn{1}{r|}{} & \multicolumn{1}{r|}{17} & \multicolumn{1}{r|}{901} & \multicolumn{1}{r|}{57.893} & \multicolumn{1}{r|}{21621.009} \\
		\cline{7-10}          &       &       &       &       &       &       &       &       &  \\
		&       &       &       &       &       &       &       & b)    &  \\
	\end{tabular}%
	\label{tab:addlabel}%
\end{table}%

%\input{FileTest_40.tex}
% Table generated by Excel2LaTeX from sheet 'FileTestHT40'
\begin{table}[htbp]
	\centering
	\caption{The test results for \textit{n} = 40 with the 2nd type constraint}
	\begin{tabular}{rrrrrrrrcr}\\
		\cline{1-4}\cline{7-10}    \multicolumn{1}{|c|}{No.} & \multicolumn{1}{c|}{Step} & \multicolumn{1}{c|}{Time} & \multicolumn{1}{c|}{roA} &       & \multicolumn{1}{c|}{} & \multicolumn{1}{c|}{No.} & \multicolumn{1}{c|}{Step} & \multicolumn{1}{c|}{Time} & \multicolumn{1}{c|}{roB} \\
		\cline{1-4}\cline{7-10}    \multicolumn{1}{|r|}{1} & \multicolumn{1}{r|}{6} & \multicolumn{1}{r|}{0.357} & \multicolumn{1}{r|}{207.869} &       & \multicolumn{1}{r|}{} & \multicolumn{1}{r|}{1} & \multicolumn{1}{r|}{4} & \multicolumn{1}{r|}{0.271} & \multicolumn{1}{r|}{31.539} \\
		\cline{1-4}\cline{7-10}    \multicolumn{1}{|r|}{2} & \multicolumn{1}{r|}{43} & \multicolumn{1}{r|}{1.078} & \multicolumn{1}{r|}{701.557} &       & \multicolumn{1}{r|}{} & \multicolumn{1}{r|}{2} & \multicolumn{1}{r|}{4} & \multicolumn{1}{r|}{0.311} & \multicolumn{1}{r|}{47.308} \\
		\cline{1-4}\cline{7-10}    \multicolumn{1}{|r|}{3} & \multicolumn{1}{r|}{69} & \multicolumn{1}{r|}{1.563} & \multicolumn{1}{r|}{1052.336} &       & \multicolumn{1}{r|}{} & \multicolumn{1}{r|}{3} & \multicolumn{1}{r|}{5} & \multicolumn{1}{r|}{0.350} & \multicolumn{1}{r|}{70.962} \\
		\cline{1-4}\cline{7-10}    \multicolumn{1}{|r|}{4} & \multicolumn{1}{r|}{107} & \multicolumn{1}{r|}{2.408} & \multicolumn{1}{r|}{1578.504} &       & \multicolumn{1}{r|}{} & \multicolumn{1}{r|}{4} & \multicolumn{1}{r|}{6} & \multicolumn{1}{r|}{0.469} & \multicolumn{1}{r|}{106.444} \\
		\cline{1-4}\cline{7-10}    \multicolumn{1}{|r|}{5} & \multicolumn{1}{r|}{163} & \multicolumn{1}{r|}{3.438} & \multicolumn{1}{r|}{2367.756} &       & \multicolumn{1}{r|}{} & \multicolumn{1}{r|}{5} & \multicolumn{1}{r|}{7} & \multicolumn{1}{r|}{0.477} & \multicolumn{1}{r|}{159.665} \\
		\cline{1-4}\cline{7-10}    \multicolumn{1}{|r|}{6} & \multicolumn{1}{r|}{373} & \multicolumn{1}{r|}{7.227} & \multicolumn{1}{r|}{5327.451} &       & \multicolumn{1}{r|}{} & \multicolumn{1}{r|}{6} & \multicolumn{1}{r|}{10} & \multicolumn{1}{r|}{0.666} & \multicolumn{1}{r|}{239.498} \\
		\cline{1-4}\cline{7-10}    \multicolumn{1}{|r|}{7} & \multicolumn{1}{r|}{561} & \multicolumn{1}{r|}{10.695} & \multicolumn{1}{r|}{7991.177} &       & \multicolumn{1}{r|}{} & \multicolumn{1}{r|}{7} & \multicolumn{1}{r|}{12} & \multicolumn{1}{r|}{0.795} & \multicolumn{1}{r|}{359.247} \\
		\cline{1-4}\cline{7-10}    \multicolumn{1}{|r|}{8} & \multicolumn{1}{r|}{843} & \multicolumn{1}{r|}{15.936} & \multicolumn{1}{r|}{11986.766} &       & \multicolumn{1}{r|}{} & \multicolumn{1}{r|}{8} & \multicolumn{1}{r|}{17} & \multicolumn{1}{r|}{1.129} & \multicolumn{1}{r|}{538.870} \\
		\cline{1-4}\cline{7-10}          &       &       &       &       & \multicolumn{1}{r|}{} & \multicolumn{1}{r|}{9} & \multicolumn{1}{r|}{23} & \multicolumn{1}{r|}{1.520} & \multicolumn{1}{r|}{808.306} \\
		\cline{7-10}          &       & \multicolumn{1}{c}{a)} &       &       & \multicolumn{1}{r|}{} & \multicolumn{1}{r|}{10} & \multicolumn{1}{r|}{47} & \multicolumn{1}{r|}{3.045} & \multicolumn{1}{r|}{1818.688} \\
		\cline{7-10}          &       &       &       &       & \multicolumn{1}{r|}{} & \multicolumn{1}{r|}{11} & \multicolumn{1}{r|}{68} & \multicolumn{1}{r|}{4.361} & \multicolumn{1}{r|}{2728.032} \\
		\cline{7-10}          &       &       &       &       & \multicolumn{1}{r|}{} & \multicolumn{1}{r|}{12} & \multicolumn{1}{r|}{100} & \multicolumn{1}{r|}{6.414} & \multicolumn{1}{r|}{4092.047} \\
		\cline{7-10}          &       &       &       &       & \multicolumn{1}{r|}{} & \multicolumn{1}{r|}{13} & \multicolumn{1}{r|}{148} & \multicolumn{1}{r|}{9.828} & \multicolumn{1}{r|}{6138.071} \\
		\cline{7-10}          &       &       &       &       & \multicolumn{1}{r|}{} & \multicolumn{1}{r|}{14} & \multicolumn{1}{r|}{220} & \multicolumn{1}{r|}{14.024} & \multicolumn{1}{r|}{9207.107} \\
		\cline{7-10}          &       &       &       &       & \multicolumn{1}{r|}{} & \multicolumn{1}{r|}{15} & \multicolumn{1}{r|}{328} & \multicolumn{1}{r|}{21.236} & \multicolumn{1}{r|}{13810.660} \\
		\cline{7-10}          &       &       &       &       & \multicolumn{1}{r|}{} & \multicolumn{1}{r|}{16} & \multicolumn{1}{r|}{490} & \multicolumn{1}{r|}{31.505} & \multicolumn{1}{r|}{20715.990} \\
		\cline{7-10}          &       &       &       &       & \multicolumn{1}{r|}{} & \multicolumn{1}{r|}{17} & \multicolumn{1}{r|}{733} & \multicolumn{1}{r|}{48.019} & \multicolumn{1}{r|}{31073.984} \\
		\cline{7-10}          &       &       &       &       &       &       &       &       &  \\
		&       &       &       &       &       &       &       & b)    &  \\
	\end{tabular}%
	\label{tab:addlabel}%
\end{table}%

% Table generated by Excel2LaTeX from sheet 'FileTestTG80'
\begin{table}[htbp]
	\centering
	\caption{The test results for \textit{ n} = 80  with the 1st type constraint}
	\begin{tabular}{rrrrrrrrcr}\\
		\cline{1-4}\cline{7-10}    \multicolumn{1}{|c|}{No.} & \multicolumn{1}{c|}{Step} & \multicolumn{1}{c|}{Time} & \multicolumn{1}{c|}{roA} &       & \multicolumn{1}{c|}{} & \multicolumn{1}{c|}{No.} & \multicolumn{1}{c|}{Step} & \multicolumn{1}{c|}{Time} & \multicolumn{1}{c|}{roB} \\
		\cline{1-4}\cline{7-10}    \multicolumn{1}{|r|}{1} & \multicolumn{1}{r|}{17} & \multicolumn{1}{r|}{2.257} & \multicolumn{1}{r|}{398.858} &       & \multicolumn{1}{r|}{} & \multicolumn{1}{r|}{1} & \multicolumn{1}{r|}{6} & \multicolumn{1}{r|}{1.329} & \multicolumn{1}{r|}{46.645} \\
		\cline{1-4}\cline{7-10}    \multicolumn{1}{|r|}{2} & \multicolumn{1}{r|}{42} & \multicolumn{1}{r|}{3.590} & \multicolumn{1}{r|}{598.287} &       & \multicolumn{1}{r|}{} & \multicolumn{1}{r|}{2} & \multicolumn{1}{r|}{6} & \multicolumn{1}{r|}{1.309} & \multicolumn{1}{r|}{69.967} \\
		\cline{1-4}\cline{7-10}    \multicolumn{1}{|r|}{3} & \multicolumn{1}{r|}{80} & \multicolumn{1}{r|}{5.654} & \multicolumn{1}{r|}{897.430} &       & \multicolumn{1}{r|}{} & \multicolumn{1}{r|}{3} & \multicolumn{1}{r|}{8} & \multicolumn{1}{r|}{1.904} & \multicolumn{1}{r|}{104.951} \\
		\cline{1-4}\cline{7-10}    \multicolumn{1}{|r|}{4} & \multicolumn{1}{r|}{137} & \multicolumn{1}{r|}{8.608} & \multicolumn{1}{r|}{1346.145} &       & \multicolumn{1}{r|}{} & \multicolumn{1}{r|}{4} & \multicolumn{1}{r|}{11} & \multicolumn{1}{r|}{2.415} & \multicolumn{1}{r|}{157.426} \\
		\cline{1-4}\cline{7-10}    \multicolumn{1}{|r|}{5} & \multicolumn{1}{r|}{223} & \multicolumn{1}{r|}{12.446} & \multicolumn{1}{r|}{2019.218} &       & \multicolumn{1}{r|}{} & \multicolumn{1}{r|}{5} & \multicolumn{1}{r|}{14} & \multicolumn{1}{r|}{3.210} & \multicolumn{1}{r|}{236.139} \\
		\cline{1-4}\cline{7-10}    \multicolumn{1}{|r|}{6} & \multicolumn{1}{r|}{351} & \multicolumn{1}{r|}{18.653} & \multicolumn{1}{r|}{3028.826} &       & \multicolumn{1}{r|}{} & \multicolumn{1}{r|}{6} & \multicolumn{1}{r|}{19} & \multicolumn{1}{r|}{4.730} & \multicolumn{1}{r|}{354.208} \\
		\cline{1-4}\cline{7-10}    \multicolumn{1}{|r|}{7} & \multicolumn{1}{r|}{543} & \multicolumn{1}{r|}{29.408} & \multicolumn{1}{r|}{4543.240} &       & \multicolumn{1}{r|}{} & \multicolumn{1}{r|}{7} & \multicolumn{1}{r|}{27} & \multicolumn{1}{r|}{6.244} & \multicolumn{1}{r|}{531.312} \\
		\cline{1-4}\cline{7-10}    \multicolumn{1}{|r|}{8} & \multicolumn{1}{r|}{831} & \multicolumn{1}{r|}{43.965} & \multicolumn{1}{r|}{6814.859} &       & \multicolumn{1}{r|}{} & \multicolumn{1}{r|}{8} & \multicolumn{1}{r|}{38} & \multicolumn{1}{r|}{7.713} & \multicolumn{1}{r|}{796.969} \\
		\cline{1-4}\cline{7-10}          &       &       &       &       & \multicolumn{1}{r|}{} & \multicolumn{1}{r|}{9} & \multicolumn{1}{r|}{55} & \multicolumn{1}{r|}{11.152} & \multicolumn{1}{r|}{1195.453} \\
		\cline{7-10}          &       & \multicolumn{1}{c}{a)} &       &       & \multicolumn{1}{r|}{} & \multicolumn{1}{r|}{10} & \multicolumn{1}{r|}{80} & \multicolumn{1}{r|}{16.487} & \multicolumn{1}{r|}{1793.179} \\
		\cline{7-10}          &       &       &       &       & \multicolumn{1}{r|}{} & \multicolumn{1}{r|}{11} & \multicolumn{1}{r|}{118} & \multicolumn{1}{r|}{23.022} & \multicolumn{1}{r|}{2689.769} \\
		\cline{7-10}          &       &       &       &       & \multicolumn{1}{r|}{} & \multicolumn{1}{r|}{12} & \multicolumn{1}{r|}{175} & \multicolumn{1}{r|}{34.075} & \multicolumn{1}{r|}{4034.653} \\
		\cline{7-10}          &       &       &       &       & \multicolumn{1}{r|}{} & \multicolumn{1}{r|}{13} & \multicolumn{1}{r|}{260} & \multicolumn{1}{r|}{70.543} & \multicolumn{1}{r|}{6051.980} \\
		\cline{7-10}          &       &       &       &       & \multicolumn{1}{r|}{} & \multicolumn{1}{r|}{14} & \multicolumn{1}{r|}{388} & \multicolumn{1}{r|}{83.998} & \multicolumn{1}{r|}{9077.970} \\
		\cline{7-10}          &       &       &       &       & \multicolumn{1}{r|}{} & \multicolumn{1}{r|}{15} & \multicolumn{1}{r|}{579} & \multicolumn{1}{r|}{108.984} & \multicolumn{1}{r|}{13616.954} \\
		\cline{7-10}          &       &       &       &       & \multicolumn{1}{r|}{} & \multicolumn{1}{r|}{16} & \multicolumn{1}{r|}{867} & \multicolumn{1}{r|}{168.209} & \multicolumn{1}{r|}{20425.431} \\
		\cline{7-10}          &       &       &       &       &       &       &       &       &  \\
		&       &       &       &       &       &       &       & b)    &  \\
	\end{tabular}%
	\label{tab:addlabel}%
\end{table}%

%\input{FileTest_80.tex}
% Table generated by Excel2LaTeX from sheet 'FileTestHT80'
\begin{table}[htbp]
	\centering
	\caption{The test results for \textit{n} = 80  with the 2nd type constraint}
	\begin{tabular}{rrrrrrrrcr}\\
		\cline{1-4}\cline{7-10}    \multicolumn{1}{|c|}{No.} & \multicolumn{1}{c|}{Step} & \multicolumn{1}{c|}{Time} & \multicolumn{1}{c|}{roA} &       & \multicolumn{1}{c|}{} & \multicolumn{1}{c|}{No.} & \multicolumn{1}{c|}{Step} & \multicolumn{1}{c|}{Time} & \multicolumn{1}{c|}{roB} \\
		\cline{1-4}\cline{7-10}    \multicolumn{1}{|r|}{1} & \multicolumn{1}{r|}{17} & \multicolumn{1}{r|}{2.424} & \multicolumn{1}{r|}{396.403} &       & \multicolumn{1}{r|}{} & \multicolumn{1}{r|}{1} & \multicolumn{1}{r|}{7} & \multicolumn{1}{r|}{1.787} & \multicolumn{1}{r|}{109.550} \\
		\cline{1-4}\cline{7-10}    \multicolumn{1}{|r|}{2} & \multicolumn{1}{r|}{43} & \multicolumn{1}{r|}{3.025} & \multicolumn{1}{r|}{594.605} &       & \multicolumn{1}{r|}{} & \multicolumn{1}{r|}{2} & \multicolumn{1}{r|}{10} & \multicolumn{1}{r|}{2.545} & \multicolumn{1}{r|}{164.325} \\
		\cline{1-4}\cline{7-10}    \multicolumn{1}{|r|}{3} & \multicolumn{1}{r|}{81} & \multicolumn{1}{r|}{4.447} & \multicolumn{1}{r|}{891.908} &       & \multicolumn{1}{r|}{} & \multicolumn{1}{r|}{3} & \multicolumn{1}{r|}{14} & \multicolumn{1}{r|}{3.285} & \multicolumn{1}{r|}{246.488} \\
		\cline{1-4}\cline{7-10}    \multicolumn{1}{|r|}{4} & \multicolumn{1}{r|}{138} & \multicolumn{1}{r|}{6.908} & \multicolumn{1}{r|}{1337.862} &       & \multicolumn{1}{r|}{} & \multicolumn{1}{r|}{4} & \multicolumn{1}{r|}{19} & \multicolumn{1}{r|}{4.677} & \multicolumn{1}{r|}{369.732} \\
		\cline{1-4}\cline{7-10}    \multicolumn{1}{|r|}{5} & \multicolumn{1}{r|}{222} & \multicolumn{1}{r|}{9.914} & \multicolumn{1}{r|}{2006.793} &       & \multicolumn{1}{r|}{} & \multicolumn{1}{r|}{5} & \multicolumn{1}{r|}{26} & \multicolumn{1}{r|}{6.597} & \multicolumn{1}{r|}{554.598} \\
		\cline{1-4}\cline{7-10}    \multicolumn{1}{|r|}{6} & \multicolumn{1}{r|}{348} & \multicolumn{1}{r|}{15.201} & \multicolumn{1}{r|}{3010.189} &       & \multicolumn{1}{r|}{} & \multicolumn{1}{r|}{6} & \multicolumn{1}{r|}{38} & \multicolumn{1}{r|}{8.805} & \multicolumn{1}{r|}{831.898} \\
		\cline{1-4}\cline{7-10}    \multicolumn{1}{|r|}{7} & \multicolumn{1}{r|}{536} & \multicolumn{1}{r|}{22.813} & \multicolumn{1}{r|}{4515.283} &       & \multicolumn{1}{r|}{} & \multicolumn{1}{r|}{7} & \multicolumn{1}{r|}{56} & \multicolumn{1}{r|}{13.169} & \multicolumn{1}{r|}{1247.846} \\
		\cline{1-4}\cline{7-10}    \multicolumn{1}{|r|}{8} & \multicolumn{1}{r|}{818} & \multicolumn{1}{r|}{33.261} & \multicolumn{1}{r|}{6772.925} &       & \multicolumn{1}{r|}{} & \multicolumn{1}{r|}{8} & \multicolumn{1}{r|}{82} & \multicolumn{1}{r|}{20.989} & \multicolumn{1}{r|}{1871.770} \\
		\cline{1-4}\cline{7-10}          &       &       &       &       & \multicolumn{1}{r|}{} & \multicolumn{1}{r|}{9} & \multicolumn{1}{r|}{121} & \multicolumn{1}{r|}{26.179} & \multicolumn{1}{r|}{2807.654} \\
		\cline{7-10}          &       & \multicolumn{1}{c}{a)} &       &       & \multicolumn{1}{r|}{} & \multicolumn{1}{r|}{10} & \multicolumn{1}{r|}{179} & \multicolumn{1}{r|}{42.407} & \multicolumn{1}{r|}{4211.481} \\
		\cline{7-10}          &       &       &       &       & \multicolumn{1}{r|}{} & \multicolumn{1}{r|}{11} & \multicolumn{1}{r|}{266} & \multicolumn{1}{r|}{67.237} & \multicolumn{1}{r|}{6317.222} \\
		\cline{7-10}          &       &       &       &       & \multicolumn{1}{r|}{} & \multicolumn{1}{r|}{12} & \multicolumn{1}{r|}{398} & \multicolumn{1}{r|}{135.115} & \multicolumn{1}{r|}{9475.833} \\
		\cline{7-10}          &       &       &       &       & \multicolumn{1}{r|}{} & \multicolumn{1}{r|}{13} & \multicolumn{1}{r|}{594} & \multicolumn{1}{r|}{138.019} & \multicolumn{1}{r|}{14213.750} \\
		\cline{7-10}          &       &       &       &       & \multicolumn{1}{r|}{} & \multicolumn{1}{r|}{14} & \multicolumn{1}{r|}{890} & \multicolumn{1}{r|}{194.137} & \multicolumn{1}{r|}{21320.624} \\
		\cline{7-10}          &       &       &       &       &       &       &       &       &  \\
		&       &       &       &       &       &       &       & b)    &  \\
	\end{tabular}%
	\label{tab:addlabel}%
\end{table}%

\end{document}